\documentclass[12pt,a4paper]{article}

\usepackage[naustrian, english]{babel} 
\usepackage[T1]{fontenc}
\usepackage{lmodern}
\usepackage{textcomp}
\usepackage{amssymb} 
\usepackage{amsthm}
\usepackage{mathtools}
\usepackage[pdftex]{graphicx}
\usepackage[export]{adjustbox}
\usepackage{caption}
\usepackage{color}
\usepackage[arrow, matrix, curve]{xy}
\usepackage{lmodern}
\usepackage[utf8]{inputenc}
\usepackage{CJKutf8}
\usepackage{url}
\usepackage{authblk}
\usepackage{hyperref}
\usepackage{float}
\usepackage{listings}
\usepackage{appendix}
\usepackage{pdfpages}
\usepackage{algorithm}
\usepackage{algorithmic}
\usepackage{geometry}
\usepackage{bbm}
\usepackage{bbold}
\usepackage{graphicx}
\usepackage{wrapfig}
\usepackage{caption}
\usepackage{subcaption}
\usepackage{xargs} 
\usepackage[colorinlistoftodos,textsize=small]{todonotes}

\newcommand{\pa}{\partial}

\newcommand{\md}{\mathrm{d}}

\newtheorem{remark}{Remark}

\usepackage{tikz}
\usetikzlibrary{positioning}
\tikzset{>=stealth}

\usepackage{pdfsync}

\allowdisplaybreaks

\newcommand{\average}{{\mathchoice {\kern1ex\vcenter{\hrule height.4pt
width 6pt depth0pt} \kern-9.7pt} {\kern1ex\vcenter{\hrule height.4pt width
4.3pt depth0pt}
\kern-7pt} {} {} }}

\title{The role of beta-amyloid and tau proteins in Alzheimer's disease: a mathematical
model on graphs}

\author[1,2]{Michiel Bertsch}
\author[3]{Bruno Franchi}
\author[3]{Maria Carla Tesi}
\author[1] {Veronica Tora}
\affil[1]{\footnotesize Department of Mathematics, University of Roma ``Tor Vergata'', Roma, Italy\endgraf (bertsch@mat.uniroma2.it, veronica.tora2@unibo.it)}
\affil[2]{\footnotesize Istituto per le Applicazioni del Calcolo ``M. Picone'', Consiglio Nazionale delle Ricerche, Roma, Italy}
\affil[3]{\footnotesize Department of Mathematics, University of Bologna, Bologna, Italy\endgraf (bruno.franchi@unibo.it, mariacarla.tesi@unibo.it)}

\date{}
\begin{document}
\maketitle


\begin{abstract} 
In this Note we study a mathematical model for the progression of Alzheimer's Disease in
the human brain. The novelty of our approach consists in the representation of the brain
as two superposed graphs where toxic proteins diffuse, 
the connectivity graph which represents the neural network,
and the proximity graph  which takes into account  the extracellular space.
Toxic proteins such as $\beta$-amyloid and $\tau$ play in fact a crucial role in the development of Alzheimer's disease and,
separately, have been
targets of medical treatments. Recent biomedical literature stresses the potential impact of the synergetic action of these proteins.
We numerically test various modelling hypotheses which confirm the relevance of this synergy.

\end{abstract}

\medskip


\medskip


\section{Introduction} Alzheimer's disease (AD) is a  neurodegenerative disease characterized by  a progressive decline in memory and other cognitive functions, leading ultimately to dementia. According to the 2019 World Alzheimer Report, it is estimated that there are presently 50 million people living with AD and related disorders, and this figure is expected to increase to 150 million by 2050 due to an increasingly aged population.

AD was first described in 1907 by Alois Alzheimer who associated  AD with histopathological  hallmarks in the brain: senile plaques  and neurofibrillary tangles (NFTs). Only after 1980 it was discovered that plaques consist primarily of aggregates of amyloid beta peptides (A$\beta$)~\cite{AD2}, whereas the main constituent of neurofibrillary tangles is misfolded tau protein ($\tau$)~\cite{AD3}. 
 In 1992, Hardy and Higgins (\cite{hardy_higgins}) formulated the so-called amyloid cascade hypothesis 
for the progression of AD: ``the deposition of A$\beta$, the main component of the plaques, is the causative agent of Alzheimer's pathology and  the neurofibrillary tangles, cell loss, vascular damage, 
and dementia follow as a direct result of this deposition''. 
Subsequently, this hypothesis has been revised along the years: although                  
 senile plaques are associated with AD, their
presence is not strictly related to the severity of the disease. 
High levels of soluble A$\beta$ 
correlate better with the presence and degree of cognitive deficits. 
Indeed, diffuse amyloid plaques are commonly present in the brains of cognitively intact elderly people.
Some authors (see for instance~\cite{haass_selkoe}) overturn the traditional perspective, and
claim that large aggregates of A$\beta$ can actually be inert or even protective to
healthy neurons. Analogously, A$\beta$ monomers have been shown to lack neurotoxicity \cite{shankar_et_al} 
and have in fact been suggested to be neuroprotective \cite{zou_et_al, giuffrida_et_al}.

In addition, experimental data have shown that the amyloid cascade hypothesis
fails to provide a fully satisfactory description of the evolution of AD,
since A$\beta$ and $\tau$ seem to act in a synergistic fashion to cause cell death (see, e.g., \cite{ittner_et_al}
and \cite{ricciarelli_fedele}). 
On the basis of these results, it has been postulated that, in the AD progression, ``A$\beta$ is the trigger and $\tau$ is the bullet''
(\cite{bloom}).

Thus,
though  A$\beta$ and $\tau$ remain currently the major therapeutic targets for 
the treatment of AD (but so far  effective therapies are  lacking),  we shall see in
Section \ref{interplay} that recent literature suggests
that the interplay between the two proteins should be crucial in the development of the disease
and must be taken into account for the development of new therapies
that should not be targeted to the two proteins separately. We refer for instance
to \cite{BFRT} for a discussion on the current medical literature.

Mathematical models are the basis for computer simulations, 
the so-called \emph{in silico} research which
effectively supplements \emph{in vivo} and \emph{in vitro} research.
An exhaustive historical overview of existing mathematical models
for AD up to 2018 can be found in \cite{carbonell_et_al}. Among more recent contributions
to {\emph {macroscopic}} modeling, we mention \cite{BFMeTT},\cite{BFMPT},\cite{raj_et_al_2021}, \cite{FHL}, 
\cite{goriely_et_al_2020}, \cite{goriely_kuhl_bick}, \cite{weickenmeier_et_al}, \cite{fornari_et_al}, \cite{goriely_PLOS},
\cite{goriely_clearence}, \cite{fornari_et_al_2020}, and references therein.
Several mathematical models, their difficulties,
pros and cons are discussed in  \cite{BFRT}, where the authors
propose 
a highly flexible mathematical model aimed to take into account as many
features of the current research as possible.  For a 
detailed presentation of the connections between this model and 
 various bio-medical interpretations, we refer to \cite{BFRT}.

 In this Note we present 
a reduced form of the model in \cite{BFRT}, focused
on the study of  the interplay of A$\beta$ and $\tau$ and their toxic effect on neurons, though
with progressive completeness and complexity.
In particular, we face the mathematical difficulty that A$\beta$ is mainly found in the extracellular space
and propagates by proximity,
while the greater part of $\tau$ is located within the neurons and propagates by neural connectivity 
\cite{ahmed_et_al} (see Section \ref{parcels}).

The core of this paper is the discussion of several  numerical
simulations, which, thanks to the flexibility of the model, enable us to  compare  
different  hypotheses, such as the amyloid cascade hypothesis and
different models for the
interaction between A$\beta$ and $\tau$.

Our mathematical model relies  on a
conceptual scheme that is currently
largely accepted in bio-medical literature.
Roughly speaking, AD is
described as a sequence of the following events: (a) the first pathological event in AD onset
is the overproduction (or the inefficient clearance) of the A$\beta$ peptide.
The polymerization of the extracellular toxic
A$\beta$ leads ultimately to
plaque deposition. (b) Toxic A$\beta$ triggers  $\tau$-pathology,
inducing misfolding of the physiological $\tau$ within neurons
and its agglomeration, yielding the formation of toxic polymers and NFT.
(c) Tau pathology spreads to connected brain
regions, driving the neurodegenerative process
(we refer for instance to \cite{gallardo_holtzman}, Fig.16.1).

%
%
%
%


\section{Interplay of A$\beta$ and $\tau$ proteins in AD}\label{interplay}

Let us start by sketching the most relevant properties of the two proteins.
At the microscopic level of the neuronal membrane, monomeric A$\beta$ peptides originate from the proteolytic cleavage of a transmembrane glycoprotein, the amyloid precursor protein (APP). By unknown and partially genetic reasons, some neurons present an  imbalance between produced and cleared A$\beta$.
On the other hand, \emph{macroscopic phenomena} take place at the level of the cerebral parenchyma. The monomeric A$\beta$  diffuses through the microscopic tortuositiy of the brain and undergoes a process of agglomeration, leading eventually to the formation of  long, insoluble amyloid fibrils, which accumulate in spherical 
deposits known as senile plaques (those observed by Alzheimer in his post-mortem studies). 

The $\tau$-protein is mainly found within axons where it
stabilizes microtubules, but is also present in smaller
amounts in dendrites and in the extracellular space.  In AD and other taupathies the $\tau$ protein undergoes
two pathological transformations: hyperphosphorylation and misfolding (see e.g. \cite{goedert_spillantini2019}).
 In the disease state, the amount of hyperphosphorylated $\tau$ is at least three times higher than that in the normal brain~\cite{AD13}.  Hyperphosphorylation of $\tau$ negatively regulates the binding of $\tau$ to microtubules, compromizing microtubule stabilization and axonal transport. It also increases the capacity of $\tau$ to self-assemble and form aggregates from oligomers to fibrils, eventually leading to its deposition as NFTs  \cite{goedert_spillantini2017}. 
Furthermore, it has been  observed that excess $\tau$ aggregates can be released into the extracellular medium, to be internalized by surrounding neurons and induce the fibrillization of endogenous $\tau$; this suggests a role for $\tau$ seeding in neurodegeneration~\cite{AD14}. This is the so-called mechanisms of release and uptake 
of the $\tau$ protein.

The soluble monomers and oligomers of the two proteins move in the brain parenchyma. We  distinguish two
different mechanisms for their diffusion, one for the intracellular $\tau$ inside the connectome (by connectivity, i.e. along  neural connections in the brain), and
the other one  by proximity in the extracellular space for the extracellular A$\beta$.

In AD both A$\beta$ oligomers and misfolded $\tau$ oligomers are known to have a toxic effect on neurons
(synaptic dysfunction, neurofibrillary tangle mediated neuron loss, and behavioral deficits), though it is still not well understood
which is the precise role of each protein in the progression of the neurodegeneration. In particular, recent studies stress
that it is crucial to understand  the interplay between the two proteins (see e.g. \cite{bush_hyman}, \cite{kara_et_al}, \cite{lewis_et_al}, \cite{ittner_et_al}, \cite{pooler_et_al}).
In \cite{ittner_et_al} possible A$\beta$-$\tau$ interactions are extensively discussed, suggesting that
A$\beta$ drives $\tau$ pathology by causing hyperphosphorylation of $\tau$, which in turn mediates toxicity in neurons,
whereas
$\tau$ mediates A$\beta$ toxicity and
 A$\beta$ and $\tau$ amplify each others toxic effects.

\section{Parcellation and connectome} \label{parcels} 

The geometric setting of our mathematical model for the interactions
between A$\beta$ and $\tau$ must 
take into account the respective 
diffusion modes: by proximity (i.e. with a local character) and by connectivity
(possibly with a non-local character).
Therefore we
  identify the cerebral parenchyma  with a pair of superposed
graphs associated with a parcellation of the brain, i.e. a subdivision of the 
human cerebral cortex into a patchwork of anatomically and functionally distinct areas
(parcels).

Following the approach proposed in \cite{raj_2012}, we  consider
a parcellation $\{\Omega_i, \, i=1,\dots N\}$ of the brain and an associated network of white-matter
fiber pathways connecting these structures. As in \cite{raj_2012}, \cite{raj_et_al_2021}, \cite{goriely_kuhl_bick} we represent
this network by means	 a finite weighted graph $G:=\{V,E\}$ in which the vertices $ V=\{x_1,\dots,x_N\}$,
are identified with points $x_i\in\Omega_i$,
and represent the $i$-th cortical or subcortical gray matter structure (i.e.the $i$-th parcel), while the edges $e_{ij}\in E$
represent the connections by white-matter fiber pathways between the $i$-th structure
and the $j$-th structure. Coherently, we introduce a family of coefficients $w_{ij}^E\ge 0$  that
measure how much the $i$-th structure and the $j$-th structure are connected
and are measured by fiber tractography \cite{raj_et_al_2015}, \cite{raj_2012}.  The
coefficients $w_{ij}^E$ are said the {\sl connectivity weights} of the graph $G$
and we
call $G$ the ``connectivity graph''. 


We also consider a second graph $\Gamma:=\{V,F\}$ with the same vertices as $G$, 
by taking a new family $F$ of edges that keep into account the Riemannian distance of the vertices
and  the heterogeneity of the cerebral parenchyma. We assume that, roughly speaking, two vertices
are adjacent if they are ``close enough''. We
call $\Gamma$ the ``proximity graph'' and we associate with $\Gamma$ a family of weights
$w_{ij}^F\ge 0$  that
take into account the geodesic distance of the $i$-th structure and the $j$-th structure in the cerebral
parenchyma. The use of two superposed graphs makes it possible to consider
simultaneously the local diffusion of A$\beta$ and the non-local diffusion of $\tau$.
We refer to Remark \ref{local} for a discussion and a comparison with the model
in \cite{BFMeTT}.

The weights  $w_{ij}^E$ and $w_{ij}^F$ make possible to introduce the 
notion of (weighted) Laplacian on $G$ and $F$.
More precisely, If $x_m$ is a vertex of $V$, we set 
$$
\pi_{m}^E:= \sum_{j} w^E_{mj}>0\qquad\mbox{and}\qquad
\pi_{m}^F:= \sum_{j} w^F_{mj}>0.
$$
Associated with the graph $G$, we  can define 
 the so-called graph Laplacian operator, $\Delta_G$  as follows. Let $g(x)$ be  any function defined over the vertices of the graph. 
  Then, for any $m, j\;\text{with}\;1\leq m, j \leq N$:
  \begin{equation}\label{nv:lapl G}
  \Delta_G g(x_m)=\dfrac{1}{\pi_{m}^E}\sum_{j}( g(x_m)-g(x_j))w_{mj}^E\, .
  \end{equation}
  The graph Laplacian $  \Delta_\Gamma$ is defined analogously.


As in  \cite{raj_2012}, \cite{raj_et_al_2021},
the map of connectomes can be extracted from a dataset 
of the MRI of a cohort of healthy subjects and diffusion-weighted MRI (dMRI) scans 
acquired previously and processed with a custom pre-processing connectomics pipeline .

\section{The mathematical model}\label{model}
Arguing as in \cite{BFRT}, the mathematical model  for the interaction between A$\beta$ and $\tau$
consists of a set of equations on the graphs $\Gamma$ and $G$ for the densities of the two proteins,
describing  their aggregation, diffusion and interactions.
Moreover, to describe the evolution of the disease, we introduce a  kinetic-type equation for a function $f$ meant to describe
the health 
state of the neurons in a fixed parcel $\Omega_i$. Roughly speaking, 
$f=f(x_i,a,t)$ is the probability density of the degree of malfunctioning $a\in [0,1]$ of neurons located in the $i$-th  parcel at time $t>0$ 
and is such that $f(x_i,a,t)\,da$ represents the fraction of neurons in the $i$-th parcel  which at time $t$ have a degree of malfunctioning 
between $a$ and $a+da$. For a precise mathematical formulation in terms of probability measures, 
see~\cite{bertsch2018SIMA}. We assume that $a$ close to $0$ stands for ``the neuron is healthy'' whereas $a$ close 
to $1$ stands for ``the neuron is dead''. This parameter, although introduced for the sake of mathematical modeling, 
can be compared with medical images from Fluorodeoxyglucose PET (FDG-PET
\cite{mosconi_et_al}).


In the present Note we consider a simplified formulation of the diffusion-agglomeration
model described in  \cite{BFRT}. This
will enable us to carry on numerical simulations, enlightening the peculiar features of different
possible scenarios for the disease, as well as the role played by the mutual relationships
of the involved parameters. First of all, we ignore the presence of $\tau$ in the extracellular
space (see, e.g., \cite{AD14}) and hence the so-called mechanisms of release and uptake 
of the $\tau$ protein. These processes are carefully described in  \cite{BFRT}
with references to biomedical literature. Here, their contribution to the spreading of
the intraneuronal $\tau$ is integrated in the diffusion equation within the neuronal network.
As in \cite{BFRT}, we keep the intra-neural diffusion process of $\tau$ governed only by the Laplace operator
on the connectivity graph. Moreover, instead of considering oligomers of (almost) arbitrary
length, we divide them in 5 compartments: monomers, dimers, short oligomers, long oligomers and
plaques or tangles.

As for Smoluchowsi's system, originally, in \cite{smoluchowski} Smoluchowski introduced a system of infinite discrete differential equations (without diffusion) for the study of rapid coagulation of aerosols. Smoluchowski's theory was successively extended to cover different physical situations. In fact, this type of equations, describing the evolving densities of diffusing particles that are prone to coagulate in pairs, models various physical phenomena, such as, e.g. polymerization, aggregation of colloidal particles, formation of stars and planets as well as biological populations, behavior of fuel mixtures in engines. We refer to \cite{drake} for an  historical account. 

As far as we know, Smoluchowski's equations for the description of the agglomeration of A$\beta$ first appears in  \cite{Murphy_Pallitto} and then in \cite{AFMT} and \cite{BFMTT}.
Let us sketch qualitatively the arguments leading to these equations:
for $k\in\mathbb N$, let $P_k$ denote a polymer of length $k$, that is a set of $k$ identical particles (monomers) that is clustered but free to move collectively in a given medium. In the course of time, polymers diffuse and, if they approach each other sufficiently close, with some probability they merge into a single polymer whose size equals the sum of the sizes of the two colliding polymers. 
For sake of simplicity we admit only binary reactions. This phenomenon is called coalescence and we write formally
$$ P_k+P_j \longrightarrow P_{k+j}$$
for the coalescence of a polymer of size $k$ with a polymer of size $j$.

Smoluchowski's equations will be presented below.
On the other hand,
the progression of AD  in the $i$-th parcel is determined by the deterioration rate $v=v_i(a,t)$
through a transport equation.
 More precisely, the equation for 
$f$ has the form (see \eqref{f_source} and \eqref{v(f)} below):
\begin{equation}\label{eq:system.f}
\partial_tf+\partial_a(vf)=0 \qquad \text{in } V\times [0,\,1]\times (0,\,T],
\end{equation}
where the deterioration rate $v=v_i(a,t)$ depends on the average health state of the neurons
in $\Omega_i$ as well on the concentrations of toxic oligomers of A$\beta$
and $\tau$ in $\Omega_i$.


%
 
Thus we are lead to the following system in $V\times (0,\,T]$: if $i=1, \dots, 5$,
let us denote by $u_i(x_m,t)$ the molar concentration of A$\beta$-polymers of length $i$ at the 
$m$-th vertex of the graph at time $t$, and by  $\tau_i(x_m,t)$ the molar concentration of 
misfolded $\tau$  polymers of length $i$ at 
the $m$-th vertex of the graph at time $t$.

%
%
If $t>0$ and $x_m\in V$, then
the equation for A$\beta$ monomers is 
\begin{equation}\label{ab_monomers}
\epsilon \dfrac{\partial u_1(x_m,t)}{\partial t} =  - d_1 \Delta_{\Gamma} u_1(x_m,t)
 -\alpha u_1(x_m,t) \sum_{j=1}^{5}u_j(x_m,t) -\sigma_1 u_1(x_m,t)
 + \mathcal{F} (f)\;,
\end{equation}
where $\alpha\ge 0$ is the probability of two polymers to coalesce,
$\sigma_1\ge 0$ takes into account the clearance of monomers,
and $\mathcal F$ is a source term that will be discussed later.

If $1<i<5$, the equations for oligomers are
\begin{equation}\label{ab_oligomers}
\begin{split}
\epsilon \dfrac{\partial u_i(x_m, t)}{\partial t} &=  - d_i \Delta_{\Gamma}  u_i(x_m, t)
 +\dfrac{\alpha}{2}\sum_{j=1}^{i-1}u_j(x_m, t) u_{i-j}(x_m,t)
\\&-\alpha u_i(x_m, t)\sum_{j=1}^{5} u_j(x_m,t)-\sigma_i u_i(x_m,t),
 \end{split}
\end{equation}
where $\sigma_i\ge 0$ takes into account the clearance of oligomers.
Finally, the evolution of amyloid plaques is described by the equation
\begin{equation}\label{ab_tangles} 
\epsilon \dfrac{\partial u_5 (x_m,t)}{\partial t} = \frac{\alpha}{2}\sum_{j+k\geq 5; \  k,\,j<5}  u_j(x_m,t) u_k(x_m,t).
\end{equation}
The equation for misfolded  $\tau$-monomers  
at $t>0$ and $x_m\in V$ reads as
\begin{equation}\label{tau_monomers}
\begin{split}
\dfrac{\partial \tau_1(x_m,t)}{\partial t} =&  - d_1 \Delta_{G}\tau_1(x_m, t) 
 -\gamma \tau_1(x_m,t) \sum_{j=1}^{5}\tau_j(x_m,t) \\&+
 c s_{\tau}(x_m,t) +C_\tau \big( \sum_{i=2}^{4} u_i(x_m,t) -\bar{U}\big)^+\;,
\end{split}
\end{equation}
where, according to Braak staging (\cite{braak}) and denoting by $V_{\mathrm{seed}} $ the vertices in the entorhinal  cortex, the source term for $\tau$ is $s_{\tau}(x_m, t) = \frac{t}{\lambda}\exp{\left(-\frac{t}{\lambda}\right)}$ for  $x_m \in V_{\mathrm{seed}}$  and $s_{\tau}(x_m, t)=0$ elsewhere (see also \cite{raj_et_al_2021}). 

In addition,
$\gamma$ is the probability of two polymers to coalesce. 
If $1<i<5$, the equations for oligomers are\begin{equation}\label{tau_oligomers}
\begin{split}
 \dfrac{\partial \tau_i(x_m,t)}{\partial t} =&  - d_i \Delta_{G} \tau_i(x_m,t))
 \\&
 +\dfrac{\gamma }{2}\sum_{j=1}^{i-1}\tau_j(x_m,t) \tau_{i-j}(x_m,t)
-\gamma \tau_i(x_m,t)\sum_{j=1}^{5}\tau_j(x_m,t)
\end{split}
\end{equation}
for $1< i < 5$.
Finally, the evolution of tangles is described by the equation
\begin{equation}\label{tau_tangles} 
 \dfrac{\partial \tau_5 (x_m,t)}{\partial t} = \frac{\gamma}{2}\sum_{j+k\geq 5; \  k,\,j<5}  \tau_j(v,t) \tau_k(v,t)\;.
\end{equation}
We choose the source term $\mathcal{F}(f)$ in equation \eqref{ab_monomers} as:
\begin{equation}\label{f_monomers}
\mathcal{F}(f)=C_{\mathcal{F}}\int_0^1 (\mu_0+a) (1-a)f(x_m,a,t)\;\md a 
\end{equation}
where $f(x_m,a,t)$ is the probability density of  the degree of malfunctioning $a \in [0,1]$ of neurons located at the $m$-th cerebral region at time $t>0$.
We assume that $a$ close to $0$ stands for “the neuron is healthy” whereas $a$ close to $1$ stands for “the neuron is dead”. The function $f(x_m,a,t)$  satisfies:
\begin{equation}\label{f_source}
\begin{split}
&\pa_t f(x_m,a,t)  + \pa_a \big( v[f(x_m,a,t)] f(x_m,a,t)\big)=0\\
&f(x_m,1,t)=0 \;\quad\text{for any}\; m =1,\dots, N, t \geq 0 \\
& f(x_m,a,0)=f_{0}(x_m,a) \;\quad\text{for any}\; m =1,\dots, N,
\end{split}
\end{equation}
where 
\begin{equation}\label{v(f)}
\begin{split}
 v[f(x_m,a,t)]&= C_{\mathcal{G}}\int_0^1  (b-a)^+ f(x_m,b,t)\;\md b +C_S (1-a) \big( \sum_{i=2}^{4}u_i(x_m,t)- \bar{U}_{A\beta} \big)^+\\& + C_T (1-a)
  \big( \sum_{i=1}^{5}\tau_i(x_m,t)- \bar{U}_{\tau} \big)^+\;.
 \end{split}
\end{equation}
Here, $C_S (1-a) \big( \sum_{i=2}^{4}u_i(x_m,t)$ and $ C_T (1-a)
  \big( \sum_{i=1}^{5}\tau_i(x_m,t)- \bar{U}_{\tau} \big)^+$ are meant to represent the
  toxic effect of $A\beta$ and $\tau$, respectively.
  
  \begin{remark}\label{local}We stress that the present model differs radically from the one
  presented in \cite{BFMeTT}. Indeed, the use here of the connectivity graph
  for the diffusion of $\tau$ allows to consider non-local effects,
  while in \cite{BFMeTT} only (pseudo-)local effects are considered.
  \end{remark}

\section{Numerical simulations and discussions}\label{numerics}

First of all, in this Section we present several numerical simulations of the mathematical
model introduced in Section \ref{model}. Subsequently, we shall discuss variants
of this model corresponding to different scenarios that have been considered in the
biomedical literature.

As in  \cite{raj_2012}, \cite{raj_et_al_2021},
the map of connectomes is extracted from a dataset 
of the MRI of a cohort of healthy subjects and diffusion-weighted MRI (dMRI) scans 
acquired previously and processed with a custom pre-processing connectomics pipeline.
In our simulations, we retrieve the data from \cite{balasz_et_al1}, \cite{balasz_et_al2},
as well as the web site  https://braingraph.org.

To perform numerical simulations, the choice of parameters is  a crucial issue. 
Since the data available in the medical literature do not provide uniquely all the 
values we need, and, in addition, our model is meant to be descriptive and not predictive,
we perform arbitrary but realistic choices. By ``realistic'' we mean that the outputs of
our numerical simulations are in some way comparable with clinical data. The parameters
that will be kept fixed are summarized in Table \ref{table:1}.

\begin{table}[h!]
	\centering
	\caption{Values of the fixed parameters}
	\begin{tabular}{c c c c c c c c c c c c c c c c} 
		  $N$ & $d_i$ & $\sigma_i$ & $\epsilon$ &  $\gamma $   & $\lambda$  & $\bar{U}$ & $C_{\mathcal{G}}$  & $C_S$ &$C_T$    & $\bar{U}_{A\beta}$ & $\bar{U}_{\tau}$ &$C_{\mathcal{F}}$ & $\mu_0$\\ [1ex] 
		\hline
		1015 &1/i &1/i & 0.1  & 4  & 10 & 0.001 & 0.1 &0.01 & 0.01 & 0.001 & 0.001 &10 & 0.01 \\  [1ex]
		
		\hline
	\end{tabular}
	
	\label{table:1}
\end{table}
The remaining parameters will be discussed later, testing various hypotheses concerning the progression of neurodegenerative processes
associated with different choices of the  constants $\alpha, C_{\tau}, c$. 
These parameters control respectively the A$\beta$ agglomeration, the production of monomeric 
$\tau$ driven by A$\beta$ oligomers and the $\tau$ seeding at the entorhinal cortex (modeled through the
Gamma-shaped functions  $s_{\tau}$).

As for the initial data, we choose  $f_0(x_m,a) $ of the form
of a function approximating the Dirac delta centered in a point $a$ close
to $0$. Thus $f_0$ represents an (almost) healthy brain.
Moreover, again in the spirit of starting from a healthy brain,
we choose
\begin{equation}\label{tau_ab_0}
\begin{split}
& u_1(x_m,0)= u_{0,1}(x_m)<<1 \; \text{for any}\; x_m \in V\\
&u_{i}(x_m,0)=0\; \text{for any}\; x_m \in V, \;2\leq i \leq 5\\
&\tau_i(x_m,0)=0 \; \text{for any}  \;x_m \in V, \;1\leq i \leq 5\;.
\end{split}
\end{equation}


\medskip

In order to provide global pictures of the evolution of the disease 
in the brain (or in large portions of the brain), we introduce some
macroscopic quantities (the total burden of A$\beta$ and $\tau$
and the average degree of malfunctioning):
assuming that all parcels have the same volume, the global amounts of A$\beta$ and $\tau$ polymers 
are given (up to  dimensional constants) by
\begin{equation}\label{global_ab_tau}
\begin{split}
&u_i(t) =\frac{1}{N}\sum_{x_m \in V} u_{i}(x_m,t)\; \text{for}\;1\leq i \leq 5\\
&\tau_i(t)=\frac{1}{N}\sum_{x_m \in V}\tau_i(x_m,t)\; \text{for}  \;1\leq i \leq 5\;,
\end{split}
\end{equation}
where we recall that $N$ is the number of vertices of $G$ and $\Gamma$.
The evolution of the disease is described in each node of the brain network by:
\begin{equation}
A(x_m,t)=\int_0^1 a f(x_m,a,t) \md a  \;\text{for} \;x_m \in V\;.
\end{equation}
Thus, the evolution of the disease in the whole brain in given by:
\begin{equation}
\label{A}
A(t)=\frac{1}{N}\sum_{x_m \in V} A(x_m,t),
\end{equation}

\vfill\newpage

\subsection{Total burden of A$\beta$ and $\tau$ in the whole brain}
\label{total burden}

\begin{figure}[h]
\centering 
\includegraphics[width=0.9\linewidth,height =0.6 \textheight]{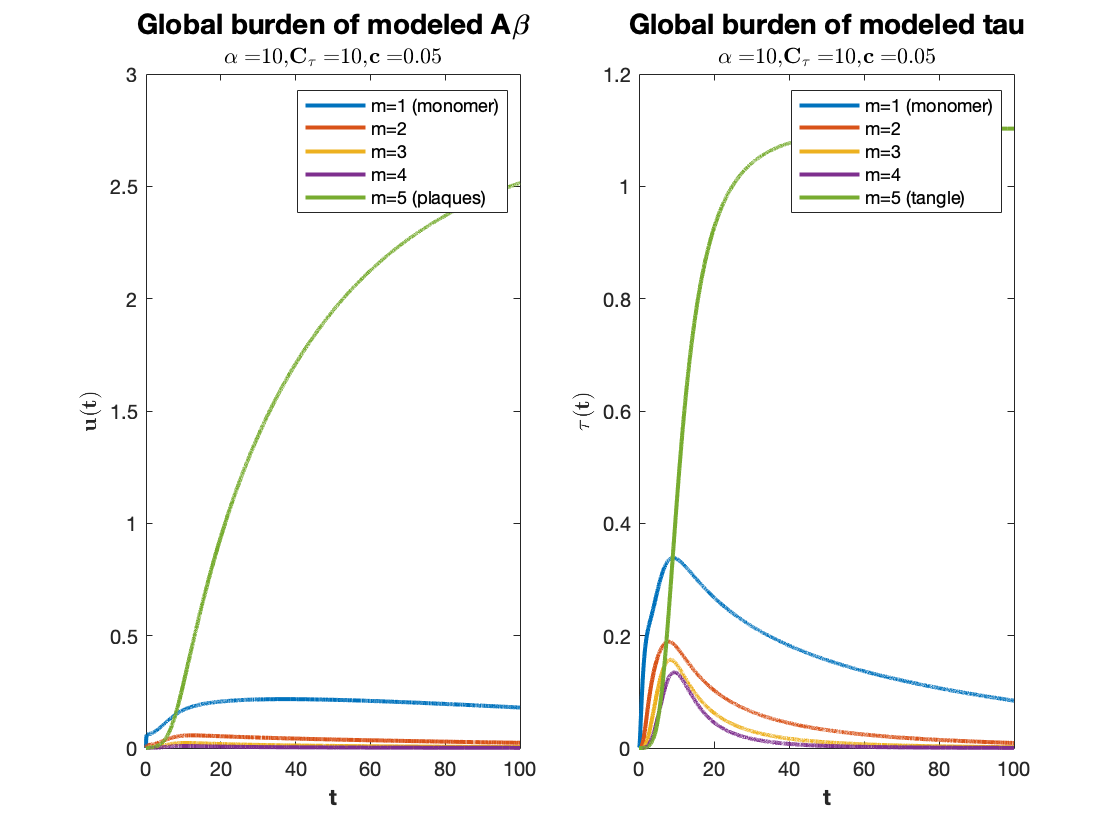} \\
\caption{\small{Temporal evolution of the global amount of  A$	\beta$ polymers (left) and $\tau$ polymers (right) of length i=1,2,3,4,5 normalized over all the vertices of the network, where we choose  $\alpha=10, C_\tau=10, c=0.05$ }.}
	\label{fig:ab_tau_caseC}
\end{figure}

In figure \ref{fig:ab_tau_caseC}, we plot the longitudinal graphs 
of $u_i(t)$ and $\tau_i(t)$ (see \eqref{global_ab_tau}) with
$\alpha=10$, $C_\tau=10$, and $c=0.05$

For both proteins, monomers' curves are the first to rise followed by those of the oligomers in increasing length. Each monomeric and oligomeric curve peaks and subsequently begins to decrease. This corresponds
to the clinical experience of advanced AD (see, e.g. \cite{ballard_et_al}).
 Moreover, low concentration of
 A$\beta$ in CSF (Celebral Spinal Fluid) is listed among diagnostic criteria and differential
diagnosis of Alzheimer's disease from other dementias.
This behavior is consistent with several factors. First, as for the A$\beta$, the imbalance between the source $\mathcal{F}(f)$ and the clearance leads monomers to increase, while concerning $\tau$-protein,  the monomers' growth is due to the seeding at EC and the effect of A$\beta$ oligomers. The first process drives the production of misfolded tau monomers at earlier times but then
  declines compatibly with the fact that the available pool of cleavable protein is limited, due to the loss of neurons. On the other hand, the production of misfolded $\tau$-monomers, governed mainly by the effect of by A$\beta$ oligomers   is lagged in time  with respect to Gamma-shaped seeding since oligomers  have to reach a certain amount to  damage the neuron. 
 The coagulation process described by means of Smoluchowski equation causes progressively the formation of larger oligomers leading to insoluble aggregates, taking active oligomers out of circulation. Insoluble clusters (plaques and tangles) are the last to develop and their curves show an increase while oligomers decrease.

\vfill\newpage

\subsection{Local burden of A$\beta$ and $\tau$ in different regions}
\label{local burden}

The aim of this Section is to simulate the progress of the formation of 
oligomers, plaques
and tangles in different regions of the cerebral parenchyma. To this end, we localize
the total burden functions of \eqref{global_ab_tau} into different regions as follows.
Let $\tilde{R}=\{ R_1, \dots, R_{\ell} \}$ be a finite partition of $V$, where each subset $R_{j}$ can be considered as a region of the brain.
Thus, for the $j$-th region, setting  $r_j=|R_{j}|$, we define:
\begin{equation}\label{local_ab_tau_R}
\begin{split}
&u_{i,R_j}(t) =\frac{1}{r_j}\sum_{x_m \in R_j} u_{i}(x_m,t)\; \text{for}\;1\leq i \leq 5\\
&\tau_{i,R_j}(t)=\frac{1}{r_j}\sum_{x_m \in R_j}\tau_i(x_m,t)\; \text{for}  \;1\leq i \leq 5\;.
\end{split}
\end{equation}
\begin{equation}\label{AR}
A_{R_j}(t)=\frac{1}{r_j}\sum_{x_m \in R_{j}} A(x_m,t).
\end{equation}

\begin{figure}[h]
	\centering
	\includegraphics[width=0.8\linewidth, height =0.45 \textheight] {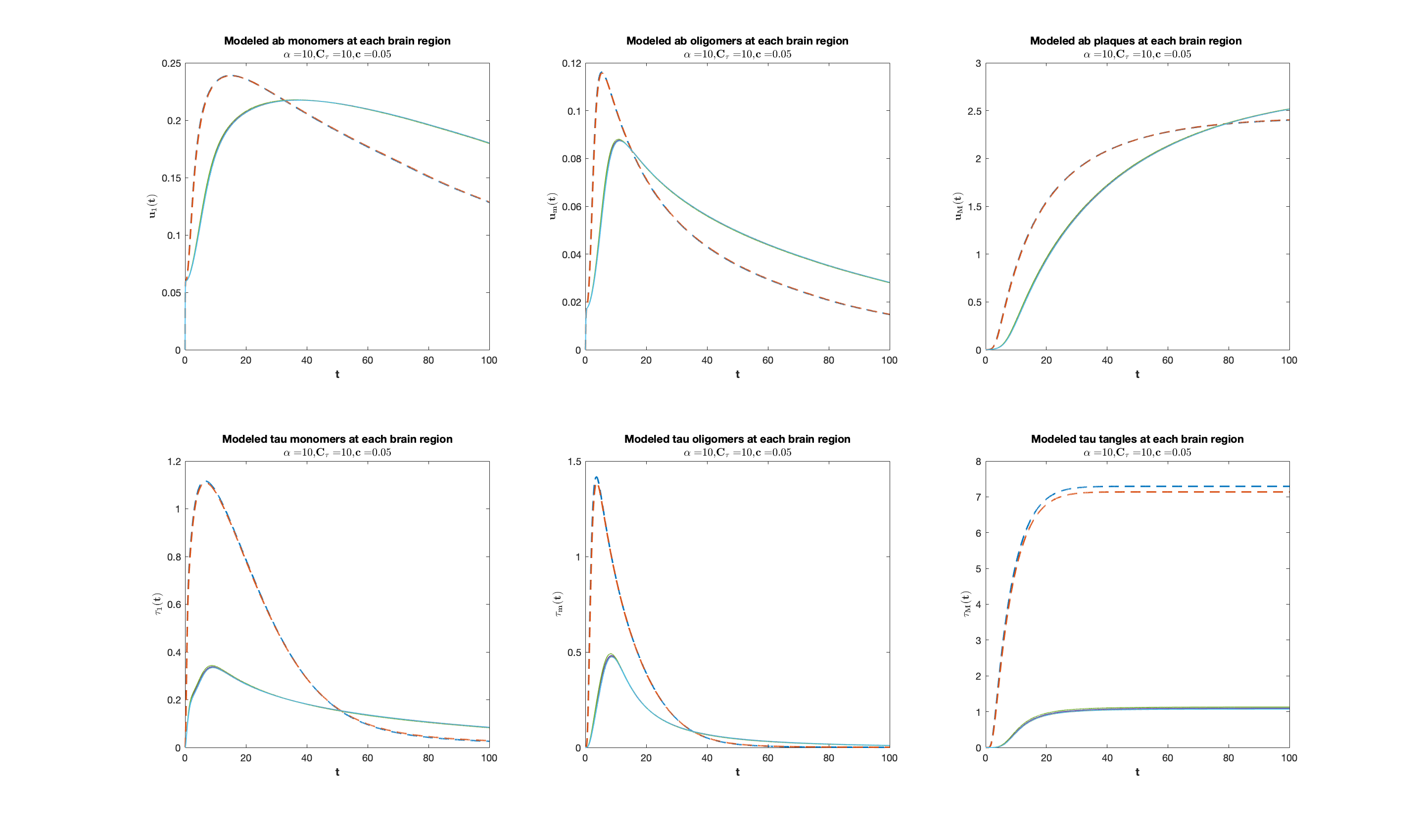}\\
\caption{\small{Temporal evolution of A$\beta$  (top) and $\tau$ (bottom)  monomers (left), oligomers (middle) and tangles (right) at each brain region for  $\alpha=10, C_\tau=10, c=0.05$. Each curve corresponds to a distinct region. The dashed curves represent the right entorhinal cortex (blue) and the left entorhinal  cortex (red).   } }
	\label{fig:regional_ab_tau_caseC}
\end{figure}

The evolution in time of monomeric, oligomeric and insoluble clusters of A$\beta$ and 
$\tau$-proteins in each brain region are shown in Figure \ref{fig:regional_ab_tau_caseC} 
for
$\alpha=10$, $C_\tau=10$, and $c=0.05$. Each curve corresponds to a distinct region. The dashed lines indicate the regions of the entorhinal cortices (EC), more precisely the right entorhinal cortex and the left entorhinal cortex.  Consistently with the plots in figure \ref{fig:ab_tau_caseC},  monomers and oligomers's regional curves rise,  peak and subsequently begin to decline while plaques and tangles' regional curves increase until they reach a plateau. 
Longitudinal graphs for different brain's regions of both A$\beta$ and $\tau$-proteins overlap, except for the entorhinal  cortex. Here, as expected, the concentrations of both soluble A$\beta$ and soluble $\tau$ are bigger than in other brain areas at  earlier times. Indeed,  in the  EC monomeric  $\tau$ and A$\beta$  rise from zero faster and consequently oligomers as well as insoluble clusters form earlier than in other brain areas. While the plaques burden in EC is bigger at least at earlier times, the tangle burden  in EC is more sizable than in other areas during all the course of the disease. Although   regional differences in the distribution of A$\beta$ and $\tau$ distribution as well as in the disease's evolution are not very evident except than in the seeding zone, small variations in A$\beta$ and $\tau$ burdens   are present. Next, the  regions near the entorhinal cortex where the burden of A$\beta$ and $\tau$ clusters is more considerable,  are (in order): amygdala, hippocampus, temporal pole, isthmus of cingulate cortex, insula cortex and parahippocampal cortex). Consistently, in these regions the evolution of the disease is more severe.

\vfill\newpage

\subsection{Global and regional behavior under different medical hypotheses}

We then test various hypotheses concerning the progression of neurodegenerative processes
which are associated to different choices of the  constants $\alpha, C_{\tau}, c$. 
We recall that these parameters control respectively the A$\beta$ agglomeration, the production of monomeric 
$\tau$ driven by A$\beta$ oligomers and the $\tau$ seeding at the entorhinal cortex.

\begin{table}[h!]
	   	\centering
		\caption{Simulated cases and respective values of parameters}
	\begin{tabular}{|c c c c|} 
		\hline
		  & $\alpha$ & $C_\tau$ & $c$\\ [1ex] 
		\hline\hline
		\text{Case A} & 10 & 0 & 0 \\ 
		\text{Case B} & 10 & 0 & 0.05\\
	    \text{Case C} & 10 & 10 & 0.05 \\
		\text{Case D} & 10 & 10 & 0 \\ 
		\text{Case E} & 0 & 10 & 0.05 \\ [1ex]
		\hline
	\end{tabular}

	\label{table:2}
\end{table}

The situations considered are listed in table \ref{table:2}. In case A the sources of  $\tau$-monomers are taken out; hence, $\tau$-protein is not produced. This situation might represent a simplified amyloid cascade hypothesis. In cases B and D, we
remove, one by one, the sources of $\tau$-monomers, represented respectively by $\beta$-amyloid oligomers and by $\tau$ seeding at EC, while in case C all these processes are included. Notice that case C has already been considered in Sections \ref{total burden}
and \ref{local burden}. In case E,  we inhibit the agglomeration of A$\beta$ in clusters and consequently,  the process of $\tau$ misfolding driven by A$\beta$ oligomers is not triggered independently from the value of the constant $C_\tau$.  This might represent a possible effect of a drug.  

In figure  \ref{fig:AD_loc_evolution}, longitudinal graphs of $A_R(t)$ are shown in each of  cases summarized in table \ref{table:2}. The curves increase in time in each brain region, consistently with the fact that the cerebral damage grows as long as the AD pathology progresses. Each curve indicate a distinct region. The dashed lines correspond to the regions of the entorhinal cortices (more precisely the right entorhinal cortex and the left entorhinal cortex) where $\tau$ misfolding is seeded.
In case A, the neuronal damage is mild and  roughly homogeneous in all regions, thus, the disease progression is slow almost everywhere. 
In  case B, the  injury is not critical in most of the cerebral areas, but, compared to case A, it is much more serious at the areas of the entorhinal cortex,  where monomeric $\tau$ is seeded. Compared to cases A and B, in cases C and D the evolution of the disease is more severe in all cerebral regions, and, in particular, in case C the EC is the most damaged brain area. Case E resembles case B, since the harm is moderate in all regions, except in the EC.   
 In such regions, the AD  pathology starts and its progression in time is more severe than in other areas of the brain.
 
 It is interesting to observe that  when the production of monomeric $\tau$ driven by oligomeric A$\beta$ is  neglected  (cases B and E), the neuronal damage is mild and  localized to some brain areas. As expected, it is serious in the region of entorhinal  cortex. On the other hand,  when the oligomeric A$\beta$ directly induces the  production of monomeric $\tau$ (case C and D), the evolution of the disease is severe in all the cerebral regions.

\begin{figure}
	\centering
	\includegraphics[width=0.9\linewidth] {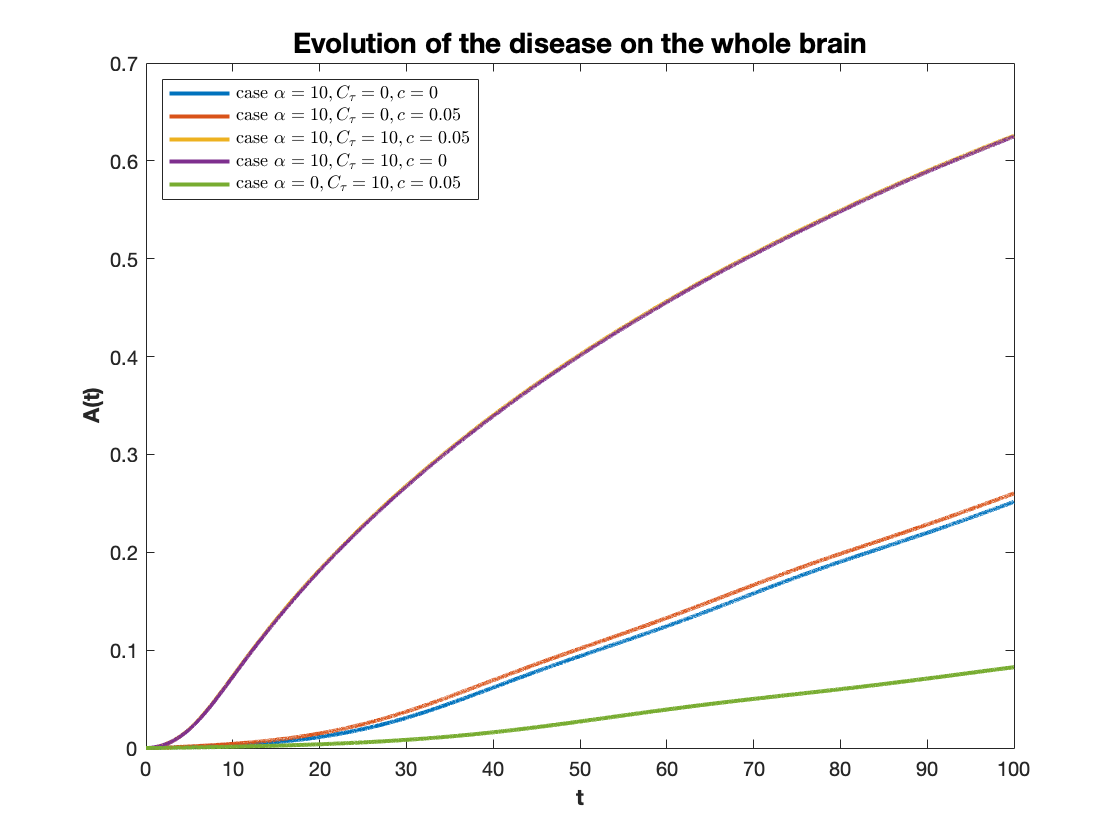}
	\caption{Temporal evolution of the AD in the whole brain  for  $\alpha=10, C_\tau=0, c=0$, (blue curve),  $\alpha=10, C_\tau=0, c=0.05$, (red curve),  $\alpha=10, C_\tau=10, c=0.05$, (purple curve), $\alpha=10, C_\tau=10, c=0$, (yellow curve), $\alpha=0, C_\tau=0, c=0.05$, (green curve),}
	\label{fig:aallcases}
\end{figure}

\begin{figure}
	\centering
	
	\text{\textbf{\tiny{ Case A}}}\\
	\includegraphics[width=0.4\linewidth, height =0.15 \textheight]{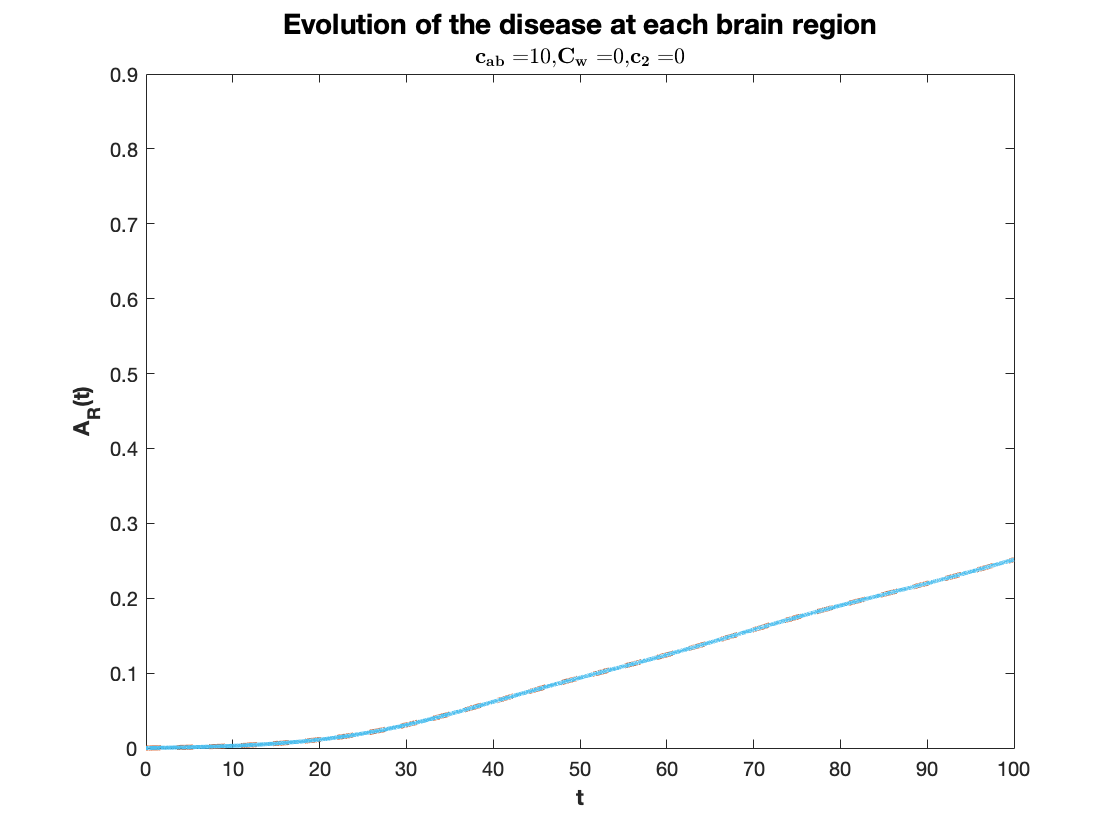}\\
	\text{\textbf{\tiny{Case B}}}\\
	\includegraphics[width=0.4\linewidth, height =0.15 \textheight]{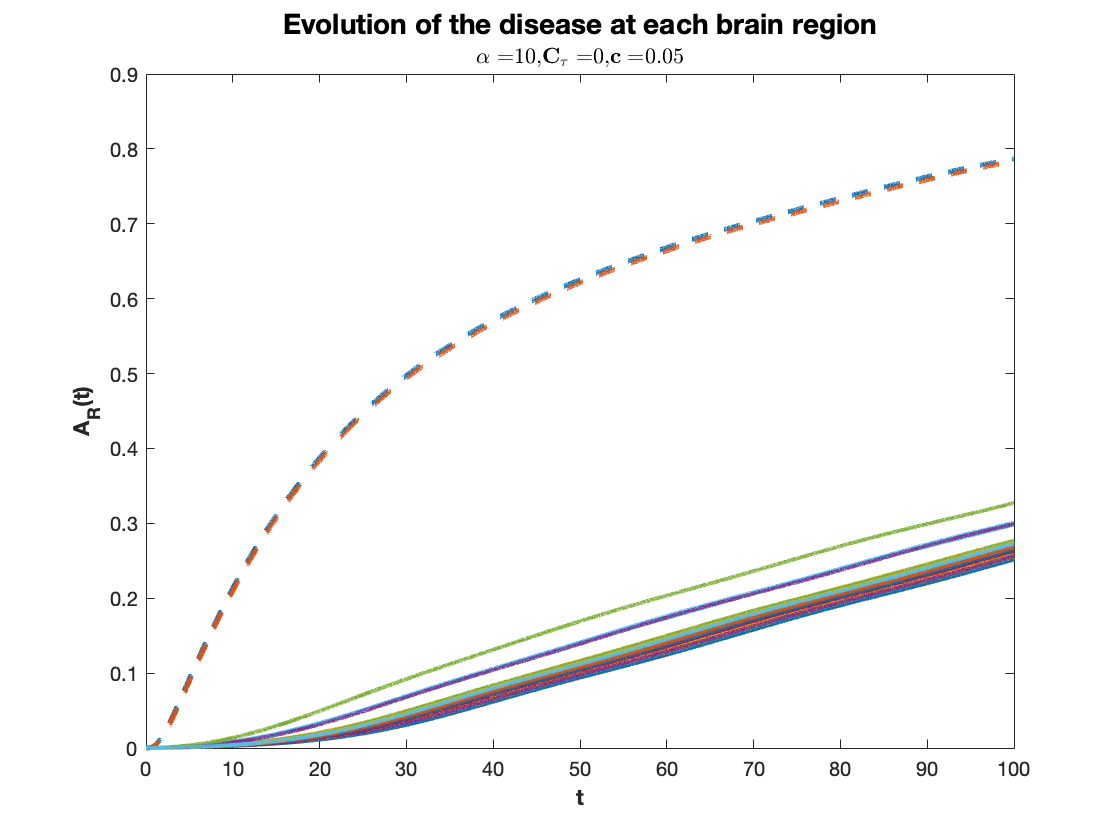}\\
	\text{\textbf{\tiny{Case C}}}\\
	\includegraphics[width=0.4\linewidth, height =0.15 \textheight]{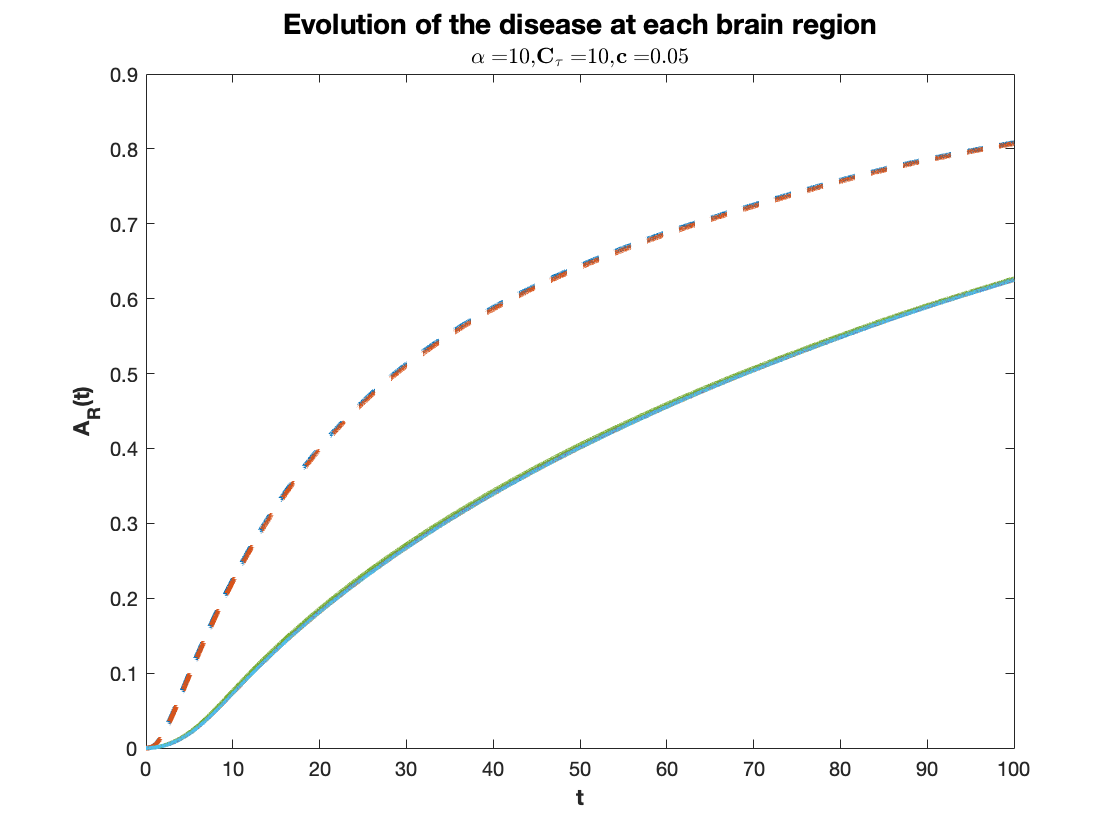}\\
	\text{\textbf{\tiny{Case D}}}\\
	\includegraphics[width=0.4\linewidth, height =0.15 \textheight]{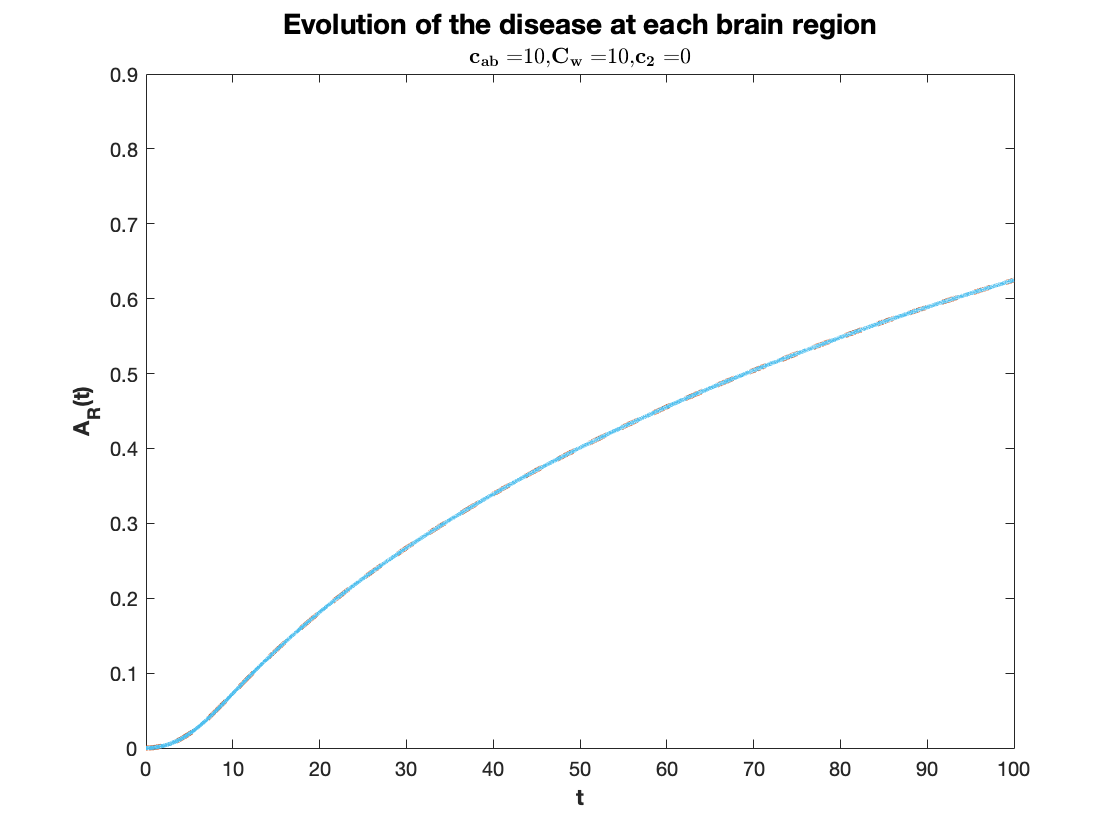}\\
	\text{\textbf{\tiny{Case E}}}\\
	\includegraphics[width=0.4\linewidth, height =0.15 \textheight]{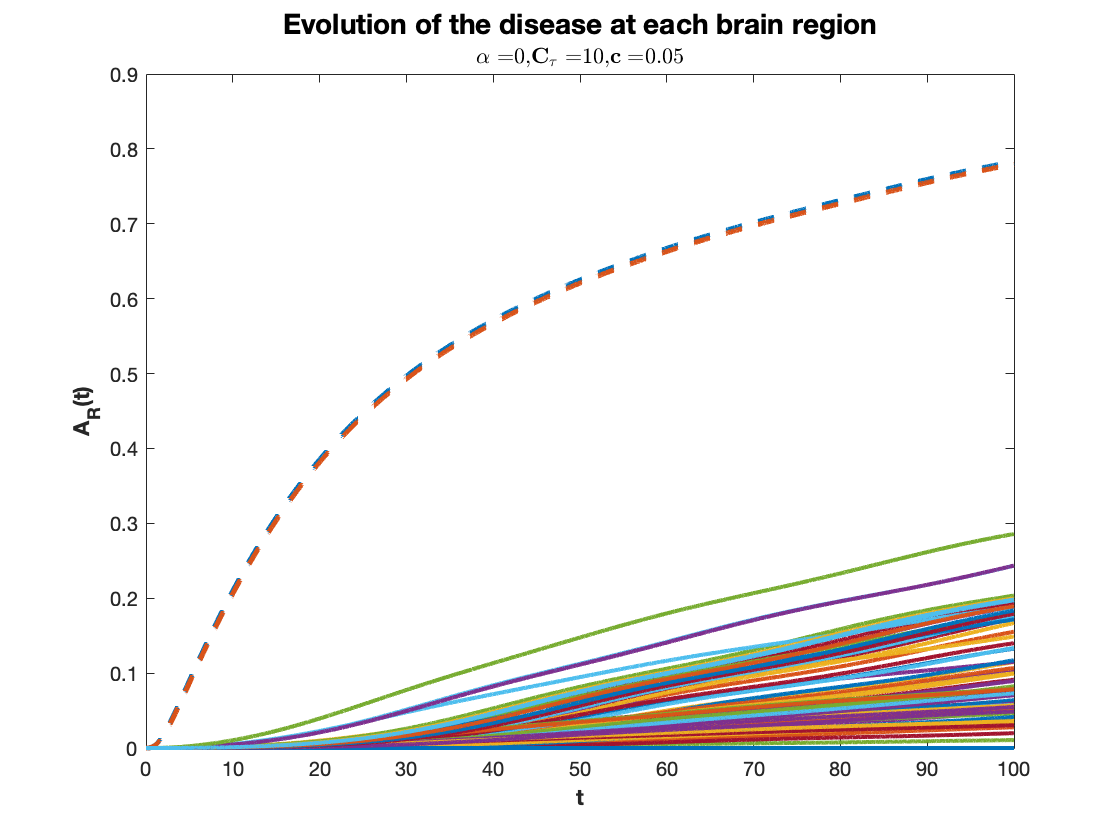}
	\caption {\small{Temporal evolution of the AD at each brain region for  $\alpha=10, C_\tau=0, c=0$, (\textbf{A}),  $\alpha=10, C_\tau=0, c=0.05$, (\textbf{B}),  $\alpha=10, C_\tau=10, c=0.05$, (\textbf{C}), $\alpha=10, C_\tau=10, c=0$, (\textbf{D}), $\alpha=0, C_\tau=0, c=0.05$, (\textbf{E}),. Each curve corresponds to a distinct region. The dashed curves represent the right entorhinal cortex (blue) and the left entorhinal  cortex (red). }}
	\label{fig:AD_loc_evolution}
\end{figure}

  In  cases B, C, E the $\tau$-seeding region (entorhinal  cortex) remains the most damaged area, followed by  amygdala, hippocampus, temporal pole, isthmuscingulate cortex, insula cortex and parahippocampal cortex  although the  differences in the evolution of AD  between regions are minimal (excluded the seeding zone) in case C and more relevant in cases B and E.
 In case D, tau seeding at EC is not considered;  consistently, the evolution of the disease at EC is not very different from any other cerebral region, suggesting that the effect of this process targets specific brain areas, that are the EC and a few more regions connected with it like amygdala and hippocampus. This fact is more evident in  cases B and E, where the tau seeding at EC is the only source for monomeric tau and in addition to EC, the most impaired cerebral zones are those connected to EC, as amygdala and hippocampus, and a few more areas connected with the latter, like, temporal pole, isthmus cingulate cortex, insula cortex and parahippocampal cortex.
In figure \ref{fig:aallcases} longitudinal graphs of $A(t)$, describing the evolution of the disease in the whole brain, are shown for each of the cases listed in table \ref{table:2}. Within the cases considered, the neuronal damage is the most severe in case C (yellow curve), followed with minimal differences by case D (purple curve), indicating, consistently, that the effect of the gamma-shaped seeding affects certain cerebral regions rather than the whole cerebral network. In  cases A (blue curve) and B (red curve), the overall cerebral damage is milder than in  cases C and D. This fact is 
interesting since it expresses the phenomenon that the injury exerted by the A$\beta$ or by both A$\beta$ and $\tau$ when the first does not directly induce the production of the latter, is much less severe than in the case  of misfolded tau injection enhanced by A$\beta$.  Thus, these observations suggest that A$\beta$ and $\tau$ are more harmful when they act together.
The simulation in figure \ref{fig:aallcases}  may provide support to the fact that the interaction between  $\tau$ and A$\beta$ increases the toxicity of both proteins and, combined  with the process of $\tau$ seeding,  plays a crucial role in brain damaging. 

In  case E, when monomeric  A$\beta$ does not coagulate in longer clusters, the cerebral damage is less serious than in the other cases, confirming once again the central role of A$\beta$ oligomers in neuronal impairment also through their involvement in $\tau$ misfolding.  Case E may represent the theoretical action of a drug able to completely inhibit the aggregation of A$\beta$. Although such action seems to be successful in limiting the overall brain's damage, it is ineffective in controlling  the impairment of specific regions (see figure  \ref{fig:AD_loc_evolution} case E ) where the disease becomes more and more serious as time goes on. 

\vfill\newpage


%
%
We observe that in case E summarized in table \ref{table:2} we inhibit the agglomeration of A$\beta$ in clusters and consequently,  the process of $\tau$ seeding driven by A$\beta$ oligomers is not triggered independently of the value of the constant $C_\tau$.  

\vfill\newpage

%

\end{document}